\documentclass[11pt]{article}
\usepackage[T1]{fontenc}
\usepackage[french]{babel}

\newcommand{\eps}{\varepsilon}
\newcommand{\var}{{\rm var}}
\newcommand{\R}{I\!\!R}
\renewcommand{\l}{{\cal L}}
\newcommand{\h}{{\cal H}}
\newcommand{\kt}{{\tilde K}}
\newcommand{\cv}{{\cal V}}
\newcommand{\tb}{{\overline t}_n}
\newcommand{\fb}{\overline {\Phi}}
\newcommand{\ov}{{\overline {\cal V}}}
\newcommand{\lp}{\stackrel{\rm P}{\to}}

\begin{document}

\centerline {\bf KERNEL ESTIMATION OF DENSITY LEVEL SETS}

\bigskip
\bigskip
\centerline {{Benoît CADRE}{\footnote{cadre@math.univ-montp2.fr}}}

\bigskip
\centerline {Laboratoire de Mathématiques, Université Montpellier II,}

\centerline {CC 051, Place E. Bataillon, 34095 Montpellier cedex 5, FRANCE}

\bigskip
\bigskip
\noindent {\bf Abstract.} Let $f$ be a multivariate density and $f_n$ be a kernel estimate of $f$ drawn from the $n$-sample $X_1,\cdots,X_n$ of i.i.d. random variables with density $f$. We compute the asymptotic rate of convergence towards 0 of the volume of the symmetric difference between the $t$-level set $\{f\geq t\}$ and its plug-in estimator $\{f_n\geq t\}$. As a corollary, we obtain the exact rate of convergence of a plug-in type estimate of the density level set corresponding to a fixed probability for the law induced by $f$.

\bigskip
\noindent {\bf Key-words:} Kernel estimate, Density level sets, Hausdorff measure. 

\bigskip
\noindent {\bf 2000 Mathematics Subject Classification:} 62H12, 62H30.

\bigskip
\noindent {\bf 1. Introduction.} Recent years have witnessed an increasing  interest in estimation of  density level sets and in related multivariate mappings problems. The main reason is the recent advent of powerfull mathematical tools and computational machinery that render these problems much more tractable. One of the most powerful application of density level sets estimation is in unsupervised {\it cluster analysis} (see Hartigan [1]), where one tries to break a complex data set into a series of piecewise similar groups or structures, each of which  may then be regarded as a separate class of data, thus reducing overall data compexity. But there are many other fields where the knowledge of density level sets is of great interest. For example, Devroye and Wise [2], Grenander [3], Cuevas [4] and Cuevas and Fraiman [5] used density support estimation for pattern recognition and for detection of the abnormal behavior of a system. 

\bigskip

In this paper, we consider the problem of estimating the $t$-level set ${\cal L}(t)$ of a multivariate probability density $f$ with support in $\R^k$ from independent random variables $X_1,\cdots,X_n$ with density $f$. Recall that for $t\geq 0$,  the $t$-level set of the density $f$ is defined as follows :
$${\cal L}(t)=\{x\in\R^k\, :\, f(x)\geq t\}.$$

\bigskip

The question now is how to define the estimates of $\l(t)$ from the $n$-sample $X_1,\cdots,X_n$ ? Even in a nonparametric framework, there are many possible answers to this question, depending on the restrictions one can impose on the level set and the density under study. Mainly, there are two families of such estimators : the {\it plug-in} estimators and the estimators constructed by an {\it excess mass} approach. Assume that an estimator $f_n$ of the density $f$ is available. Then a straightforward estimator of the level set $\l(t)$ is $\{f_n\geq t\}$, the plug-in estimator. Molchanov [6, 7] and  Cuevas and Fraiman [5] proved consistency of these estimators and obtained some rates of convergence. The excess mass approach suggest to first consider the empirical mapping $M_n$ defined for every borel set $L\subset\R^k$ by 
$$M_n(L)=\frac{1}{n} \sum_{i=1}^n {\bf 1}_{\{X_i\in L\}}-t \lambda(L),$$
where $\lambda$ denotes the Lebesgue measure on $\R^k$. A natural estimator of $\l(t)$ is a maximizer of $M_n(L)$ over a given class of borel sets $L$. For different classes of level sets (mainly star-shaped or convex level sets), estimators based on the excess mass approach were studied by Hartigan [8], M\"uller [9], M\"uller and Sawitzki [10], Nolan [11] and  Polonik [12], who proved consistency and found certain rates of convergence. When the level set is star-shaped, Tsybakov [13] recently proved that the excess mass approach gives estimators with optimal rates of convergence in an asymptotically minimax sense, whithin the studied classes of densities. Though this result has a great theoretical interest, assuming the level set to be convex or star-shaped appears to be somewhat unsatisfactory for the statistical applications. Indeed, such an assumption does not permit to consider the important case where the density under study is multimodal with a finite number of modes, and hence the results can not be applied to cluster analysis in particular.  In comparison, the plug-in estimators do not care about the specific shape of the level set. Moreover, another advantage of the plug-in approach is that it leads to easily computable estimators. We emphasize that, if the excess mass approach often gives estimators with optimal rates of convergence, the complexity of the computational algorithm of such an estimator is high, due to the presence of the maximizing step (see the computational algorithm proposed by Hartigan, [8]). 

\bigskip
In this paper, we study a plug-in type estimator of the density level set $\l(t)$, using a kernel density estimate of $f$ (Rosenblatt, [14]). Given a kernel $K$ on $\R^k$ ({\it i.e.}, a probability density on $\R^k$) and a bandwidth $h=h(n) >0$ such that $h\to 0$ as $n$ grows to infinity, the kernel estimate of $f$ is given by 
$$f_n(x)=\frac{1}{nh^k}\sum_{i=1}^n K\Big(\frac{x-X_i}{h}\Big), \ x\in\R^k.$$
We let the plug-in estimate $\l_n(t)$ of $\l(t)$ be defined as
$$\l_n(t)=\{x\in\R^k\, : \, f_n(x)\geq t\}.$$

\bigskip
In the whole paper, the distance between two borel sets in $\R^k$ is a measure -in particular the volume or Lebesgue measure $\lambda$ on $\R^k$- of the symmetric difference denoted $\Delta$ ({\it i.e.}, $A\Delta B=(A\cap B^c)\cup(A^c\cap B)$ for all sets $A,B$). Our main result (Theorem 2.1) deals with the limit law of 
$$\sqrt {nh^k}\,\lambda\Big(\l_n(t)\Delta\l(t)\Big),$$
which is proved to be degenerate.

\bigskip
Consider now the following statistical problem. In cluster analysis for instance, it is of interest to estimate the density level set corresponding to a fixed probability $p\in [0,1]$ for the law induced by $f$. The data contained in this level set can then be regarded as the most important data if $p$ is far enough from 0. Since $f$ is unknown, the level $t$ of this density level set is unknown as well. The natural estimate of the target density level set $\l(t)$ becomes $\l_n(t_n)$, where $t_n$ is such that 
$$\int_{\l_n(t_n)} f_nd\lambda=p.$$
As a consequence of our main result, we obtain in Corollary 2.1  the exact asymptotic rate of convergence of $\l_n(t_n)$ to $\l(t)$. More precisely, we prove that for some $\beta_n$ which only depends on the data, one has :
 $$\beta_n\sqrt {nh^k}\,\lambda\Big(\l_n(t_n)\Delta\l(t)\Big)\to \sqrt {\frac{2}{\pi}\int K^2d\lambda}$$ 
in probability.

\bigskip
The precise formulations of Theorem 2.1 and Corollary 2.1 are given in Section 2. Section 3 is devoted to the proof of Theorem 2.1 while the proof of Corollary 2.1 is given in Section 4. The appendix is dedicated to a change of variables formula involving the $(k$-$1)$-dimensional Hausdorff measure (Proposition A).

\bigskip
\noindent {\bf 2. The main results.} 

\bigskip
\noindent {\bf 2.1 Estimation of $t$-level sets.} In the following, $\Theta\subset (0,\infty)$ denotes an open interval and $\|.\|$ stands for the euclidean norm over any finite dimensional space. Let us introduce the hypotheses on the density $f$ :
\begin{itemize}
\item [{\bf H1.}] $f$ is twice continuously differentiable and $f(x)\to 0$ as $\|x\|\to\infty$;
\item [{\bf H2.}] For all $t\in\Theta$,
$$\inf_{f^{-1}(\{t\})} \|\nabla f\|>0,$$
\end{itemize}
where, here and in the following, $\nabla\psi(x)$ denotes the gradient at $x\in\R^k$ of the differentiable function $\psi\,:\, \R^k\to\R$. Next, we introduce the assumptions on the kernel $K$ :
\begin{itemize}
\item [{\bf H3.}] $K$ is a continuously differentiable and  compactly supported function. Moreover, there exists a monotone nonincreasing function $\mu\, :\, \R_+\to\R$ such that $K(x)=\mu(\|x\|)$ for all $x\in\R^k$.
\end{itemize}
 The assumption on the support of $K$ is only provided for simplicity of the proofs. As a matter of fact, one could consider a more general class of kernels, including the gaussian kernel for instance. Moreover, as we will use Pollard's results [15], $K$ is assumed to be of the form $\mu(\|.\|)$.

\bigskip Throughout the paper, $\h$ denotes the $(k$-$1)$-dimensional Hausdorff measure on $\R^k$ ({\it cf.}  Evans and Gariepy, [16]). Recall that $\h$ agrees with ordinary ``$(k$-$1)$-dimensional surface area'' on nice sets. Moreover,  $\partial A$ is the boundary of the set $A\subset\R^k$, 
\[ \alpha(k) = \left\{ \begin{array}{ll}
3 & \mbox{if $k=1$;} \\
k+4 & \mbox{if $k\geq 2$.}
\end{array}
\right. \]
and for any bounded borel function $g\, :\, \R^k\to\R_+$, $\lambda_g$ stands for the measure defined for each borel set $A\subset\R^k$ by
$$\lambda_g(A)=\int_A g\, d\lambda.$$
Finally, the notation $\lp$ denotes the convergence in probability. 

\bigskip It can be proved that if {\bf H1}, {\bf H3} hold and if $\lambda(\partial \l(t))=0$, one has  :
$$\lambda\Big(\l_n(t)\Delta \l(t)\Big)\lp 0.$$
The aim of Theorem 2.1 below is to obtain the exact rate of convergence.

\bigskip
\noindent {\bf Theorem 2.1.} {\it  Let $g\, :\, \R^k\to\R_+$ be a bounded borel function and assume that {\bf H1}-{\bf H3} hold. If $nh^k/(\log n)^{16}\to\infty$ and  $nh^{\alpha(k)} (\log n)^2\to 0$, then for almost every (a.e.) $t\in\Theta$ :
$$\sqrt {nh^k}\, \lambda_g\Big(\l_n(t)\Delta\l(t)\Big)\lp \sqrt{\frac{2t}{\pi} \int K^2d\lambda}\, \int_{\partial \l(t)} \frac{g}{\|\nabla f\|}d\h.$$
}

\bigskip
\noindent {\bf Remarks 2.1.} $\bullet$ Notice that the rightmost integral is defined because  $g$ is bounded and $\l(t)$ is a compact set for all $t>0$ according to {\bf H1}.

\noindent $\bullet$ In practice, this result is mainly interesting when $g\equiv 1$, since we then have the asymptotic behavior of the volume of the symmetric difference between the two level sets. The general case is provided for the proof of Corollary 2.1 below. 

\noindent $\bullet$ If we only assume $f$ to be Lipschitz instead of {\bf H1}, then $f$ is an almost everywhere continuously differentiable function by Rademacher's theorem and Theorem 2.1 holds under the additional assumption on the bandwidth : $nh^{k+2}(\log n)^2\to 0$.

\bigskip
\noindent {\bf 2.2 Estimation of level sets with fixed probability.} In order to derive the corollary, we need an additional condition on $f$.
\begin{itemize}
\item[{\bf H4.}] For all $t\in (0,\sup_{\R^k} f]$,  $\lambda(f^{-1}[t-\eps,t+\eps])\to 0$ as $\eps\to 0$. Moreover, $\lambda(f^{-1}(0,\eps])\to 0$ as $\eps\to 0$.
\end{itemize}
Roughly speaking, {\bf H4} means that the sets where $f$ is constant do not charge the Lebesgue measure on $\R^k$. Many densities with a finite number of local extrema satisfy {\bf H4}. However, notice that if $f$ is a continuous density such that  $\lambda(f^{-1}(0,\eps])\to 0$ as $\eps\to 0$, then it is compactly supported.

\bigskip
Let us now denote by ${\cal P}$ the application
\[ {\cal P} : \begin{array}{ll}
[0,\sup_{\R^k} f] & \to [0,1] \\
\qquad t & \mapsto \lambda_f(\l(t)).\\
\end{array}
 \]
Observe that $\cal P$ is one-to-one if $f$ satisfies {\bf H1}, {\bf H4}. Then, for all $p\in [0,1]$, let $t^{(p)}\in [0,\sup_{\R^k} f]$ be the unique real number such that $\lambda_f(\l(t^{(p)}))=p$. Morevover, let $t_n^{(p)}\in [0,\sup_{\R^k} f_n]$ be such that $\lambda_{f_n}(\l_n(t_n^{(p)}))=p$. Notice that $t_n^{(p)}$ does exists since $f_n$ is a density on $\R^k$. 

\bigskip The aim of Corollary 2.1 below is to obtain the exact rate of convergence of $\l_n(t_n)$ to $\l(t)$. We also introduce an estimator of the unknown integral in Theorem 2.1.

\bigskip
\noindent {\bf Corollary 2.1.} {\it Let $k\geq 2$, $(\alpha_n)_n$ be a sequence of positive real numbers such that $\alpha_n\to 0$ and  assume that {\bf H1}-{\bf H4} hold. If $nh^{k+2}/\log n\to \infty$, $nh^{k+4} (\log n)^2\to 0$ and $\alpha_n^2nh^k/(\log n)^2\to \infty$ then, for a.e. $p\in{\cal P}(\Theta)$ :
$$\sqrt {nh^k}\frac {\beta_n}{\sqrt {t_n^{(p)}}}\, \lambda\Big(\l_n(t_n^{(p)})\Delta\l(t^{(p)})\Big)\lp \sqrt {\frac{2}{\pi}\int K^2d\lambda},$$
where $\beta_n=\alpha_n/\lambda\big(\l_n(t_n^{(p)})-\l_n(t_n^{(p)}+\alpha_n)\big).$
}

\bigskip
\noindent {\bf Remarks 2.2.} $\bullet$ It is of statistical interest to mention the fact that under the assumptions of the corollary, we have for all $p\in [0,1]$ : $t_n^{(p)}\to t^{(p)}$ with probability 1 (see Lemma 4.3).

\noindent $\bullet$ When $k=1$, the conditions of Theorem 2.1 on the bandwidth $h$ do not permit to derive Corollary 2.1. In practice, estimations of density level sets and their applications to cluster analysis for instance are mainly interesting in high-dimensional problems.

\bigskip
\noindent {\bf 3. Proof of Theorem 2.1.} 

\bigskip
\noindent {\bf 3.1. Auxiliary results and proof of Theorem 2.1.} For all $t>0$, let 
$$\cv_n^t=f^{-1} \Big[t-\frac{(\log n)^{\beta}}{\sqrt {nh^k}},t\Big] \quad {\rm and} \quad \ov_n^t=f^{-1} \Big[t, t+\frac{(\log n)^{\beta}}{\sqrt {nh^k}}\Big],$$
where $\beta >1/2$ is fixed. Moreover, $\kt$ stands for the real number :
$$\kt =\int K^2 d\lambda.$$

\bigskip
\noindent {\bf Proposition 3.1.} {\it Let $g\, :\, \R^k\to \R_+$ be a bounded borel function and assume that {\bf H1}-{\bf H3} hold. If $nh^k/(\log n)^{31\beta}\to\infty$ and $nh^{\alpha(k)} (\log n)^{2\beta}\to 0$, then for a.e. $t\in\Theta$ :
\begin{eqnarray*}
\lim_n \sqrt {nh^k} \int_{\cv_n^t} P(f_n(x)\geq t)d\lambda_g(x) & = & \lim_n \sqrt {nh^k} \int_{\ov_n^t} P(f_n(x)< t)d\lambda_g(x)\\
& = & \sqrt{\frac{t\kt}{2\pi}}\, \int_{\partial \l(t)} \frac{g}{\|\nabla f\|}d\h.
\end{eqnarray*}
}

\bigskip
\noindent {\bf Proposition 3.2.} {\it Let $g\, :\, \R^k\to \R_+$ be a bounded borel function and assume that {\bf H1}-{\bf H3} hold. If $nh^k/(\log n)^{5\beta}\to\infty$ and $nh^{\alpha(k)} (\log n)^{2\beta}\to 0$, then for a.e. $t\in\Theta$ :
$$\lim_n nh^k \var\,\Big[\lambda_g\Big(\cv_n^t\cap \l_n(t)\Big)\Big]=0=\lim_n nh^k \var\,\Big[ \lambda_g\Big(\ov_n^t\cap \l_n(t)^{c}\Big)\Big].$$}

\bigskip
\noindent {\bf Proof of Theorem 2.1.} Let $t\in\Theta$ be such that both conclusions of  Propositions 3.1 and 3.2 hold. According to {\bf H3} and Pollard ([15], Theorem 37 and Problem 28, Chapter II), we have almost surely (a.s.) :
$$\sup_{\R^k}|f_n-Ef_n|\to 0.$$
Moreover, since both $\sup_n Ef_n(x)$ and $f(x)$ vanish as $\|x\|\to\infty$ by {\bf H1}, {\bf H3}, we have :
$$\sup_{\R^k}|Ef_n-f|\to 0.$$
Thus, a.s.  and for $n$ large enough :
$$\sup_{\R^k} |f_n-f|\leq \frac{t}{2}.$$
Consequently, $\l_n(t)\subset \l(t/2)$ and since $\l(t)\subset\l(t/2)$, we get :
\begin{eqnarray*}
\lambda_g\Big(\l_n(t)\Delta\l(t)\Big)=\int_{\l(t/2)} {\bf 1}_{\{f_n<t,f\geq t\}}d\lambda_g+\int_{\l(t/2)} {\bf 1}_{\{f_n\geq t,f< t\}}d\lambda_g. \qquad (3.1)
\end{eqnarray*}
Let 
$$A_n=\Big\{ \sqrt{nh^k} \sup_{\l(t/2)} |f_n-f|\leq (\log n)^{\beta}\Big\}.$$
Since $\l(t/2)$ is a compact set by {\bf H1}, it is a classical exercise to prove that $P(A_n)\to 1$ under the assumptions of the theorem. Hence, one only needs to prove that the result of Theorem 2.1 holds on the event $A_n$. But on $A_n$, one has according to (3.1): $\lambda_g\big(\l_n(t)\Delta\l(t)\big)=J_n^1+J_n^2$, where :
$$J_n^1=\lambda_g\Big(\ov_n^t\cap\l_n(t)^{c}\Big) \ {\rm and} \ J_n^2=\lambda_g\Big(\cv_n^t\cap\l_n(t)\Big).$$
By Propositions 3.1 and 3.2, if $j=1$ or $j=2$ :
$$\sqrt {nh^k} J_n^j\lp \sqrt{\frac{t\kt}{2\pi}}\, \int_{\partial \l(t)} \frac{g}{\|\nabla f\|}d\h, \qquad (3.2)$$
if the bandwidth $h$ satisfies $nh^{\alpha(k)}(\log n)^{2\beta}\to 0$ and $nh^k/(\log n)^{31\beta}\to\infty$. Letting  $\beta=16/31$, the theorem is proved $\bullet$

\bigskip
\noindent {\bf 3.2. Proof of Proposition 3.1.} Let $X$ be a random variable with density $f$,
$$V_n(x)=\var\, K\Big(\frac{x-X}{h}\Big) \  {\rm and} \ Z_n(x)=\frac {h^k\sqrt n}{\sqrt {V_n(x)}}(f_n(x)-Ef_n(x)),$$
for all $x\in\R^k$ such that $V_n(x)\neq 0$. Moreover, $\Phi$ denotes the distribution function of the ${\cal N} (0,1)$ law. 

In the proofs, $c$ denotes a positive constant whose value may vary from line to line.

\bigskip
\noindent {\bf Lemma 3.1.} {\it Assume that {\bf H1}, {\bf H3} hold and let ${\cal C}\subset \R^k$ be a compact set such that $\inf_{{\cal C}}f >0$. Then, there exists $c>0$ such that for all $n\geq 1$, $x\in{\cal C}$ and $u\in\R$:
$$|P(Z_n(x)\leq u)-\Phi(u)|\leq \frac{c}{\sqrt {nh^k}}.$$
}

\bigskip
\noindent {\bf Proof.} By the Berry-Essèen inequality ({\it cf.} Feller, [17]), one has for all $n\geq 1$, $u\in\R$ and  $x\in\R^k$ such that $V_n(x)\neq 0$:
$$|P(Z_n(x)\leq u)-\Phi(u)|\leq \frac{3}{\sqrt {nV_n(x)^3}} E\Big|K\Big(\frac{x-X}{h}\Big)-EK\Big(\frac{x-X}{h}\Big)\Big|^3.$$
It is a classical exercise to deduce from {\bf H1}, {\bf H3} that 
$$\sup_{x\in{\cal C}}  E\Big|K\Big(\frac{x-X}{h}\Big)-EK\Big(\frac{x-X}{h}\Big)\Big|^3\leq c\, h^k \ {\rm and} \ \inf_{x\in{\cal C}} V_n(x)\geq c\, h^k,$$
hence the lemma $\bullet$ 

\bigskip For all borel bounded function $g\, :\, \R^k\to\R_+$, we let $\Theta_0(g)$ to be the set of $t\in\Theta$ such that :
$$\lim_{\eps\searrow 0} \frac{1}{\eps} \lambda_g\Big(f^{-1} [t-\eps,t]\Big)=
\lim_{\eps\searrow 0} \frac{1}{\eps} \lambda_g\Big(f^{-1} [t, t+\eps]\Big)=
\int_{\partial \l(t)} \frac{g}{\|\nabla f\|}d\h.$$

\bigskip
\noindent {\bf Lemma 3.2.} {\it Let $g\, :\, \R^k\to\R_+$ be a borel bounded function and assume that {\bf H1}, {\bf H2} hold. Then we have : $\Theta_0(g)=\Theta$ a.e.}

\bigskip
\noindent {\bf Proof.} According to {\bf H1}, {\bf H2}, for all $t\in\Theta$, there exists $\eta >0$ such that :
$$\inf_{f^{-1}[t-\eta,t+\eta]}\|\nabla f\|>0.$$
We deduce from Proposition A that for all $t\in\Theta$ and $\eps>0$ small enough :
$$\frac{1}{\eps} \lambda_g\Big(f^{-1} [t-\eps,t]\Big)=\frac{1}{\eps} \int_{t-\eps}^t \int_{\partial \l(s)} \frac{g}{\|\nabla f\|} d\h\, ds.$$
Using the Lebesgue-Besicovitch theorem ({\it cf.} Evans and Gariepy, [16], Theorem 1, Chapter I), we then have for a.e. $t\in\Theta$:
$$\lim_{\eps\searrow 0}\frac{1}{\eps}  \lambda_g\Big(f^{-1} [t-\eps,t]\Big)=\int_{\partial \l(t)} \frac{g}{\|\nabla f\|} d\h,$$
and the same result holds for $\lambda_g\big(f^{-1} [t,t+\eps]\big)$ instead of $\lambda_g\big(f^{-1} [t-\eps,t]\big)$, hence the lemma $\bullet$ 

\bigskip It is a straightforward consequence of Lemma 3.2 above that $\lambda(\partial \l(t))=0$ for a.e. $t\in\Theta$. For simplicity, we shall assume throughout that this is true for all $t\in\Theta$. Since $\Theta$ is an open interval, we have in particular
$$\lambda\Big(f^{-1}[t-\eps,t+\eps]\Big)=\lambda\Big(f^{-1}(t-\eps,t+\eps)\Big),$$
for all $t\in\Theta$ and $\eps >0$ small enough.

\bigskip We now let for $t\in\Theta$ and $x\in\R^k$ such that $f(x)V_n(x)\neq 0$ :
$$t_n(x)=\sqrt {\frac{nh^k}{\kt f(x)}}(t-f(x)) \ {\rm and} \ \tb(x)=\frac{h^k\sqrt n}{\sqrt {V_n(x)}}(t-Ef_n(x)),$$
and finally, $\fb(u)=1-\Phi(u)$ for all $u\in\R$. 

\bigskip
\noindent {\bf Lemma 3.3.} {\it Let $g\, :\, \R^k\to\R_+$ be a bounded borel function and assume that {\bf H1}, {\bf H2} hold. If $nh^k/(\log n)^{2\beta}\to\infty$ and $nh^{k+4} (\log n)^{2\beta}\to 0$, then for all $t\in\Theta_0(g)$ :
\begin{eqnarray*}
&  & \lim_n \sqrt {nh^k} \Big[ \int_{\cv_n^t} P(f_n(x)\geq t)d\lambda_g(x)-\int_{\cv_n^t} \fb(t_n(x))d\lambda_g(x)\Big]=0 \\
& {\rm and} & 
\lim_n \sqrt {nh^k} \Big[ \int_{\ov_n^t} P(f_n(x)< t)d\lambda_g(x) -\int_{\ov_n^t} \Phi(t_n(x))d\lambda_g(x)\Big]=0.
\end{eqnarray*}
}

\bigskip
\noindent {\bf Proof.} We only prove the first equality. Let $t\in\Theta_0(g)$. First note that for all $x\in\R^k$ such that $V_n(x)\neq 0$ :
$$P(f_n(x)\geq t)=P(Z_n(x)\geq \tb(x)).$$
There exists a compact set ${\cal C}\subset\R^k$ such that $\inf_{{\cal C}} f>0$ and $\cv_n^t\subset{\cal C}$ for all $n$. Observe that by Lemma 3.1 and the above remarks,
$$\sqrt {nh^k}\Big[ \int_{\cv_n^t} P(f_n(x)\geq t)d\lambda_g(x)-\int_{\cv_n^t} \fb(\tb(x))d\lambda_g(x)\Big]\leq c\, \lambda_g(\cv_n^t).$$
 Since $\lambda_g(\cv_n^t)\to 0$ by Lemma 3.2, one only needs now to prove that :
$$E_n:=\sqrt {nh^k} \int_{\cv_n^t}|\fb(\tb(x))-\fb(t_n(x))|d\lambda_g(x)\to 0.$$
One deduces from the Lipschitz property of $\Phi$ that
$$E_n\leq c\sqrt {nh^k} \lambda_g(\cv_n^t) \sup_{x\in\cv_n^t}|\tb(x)-t_n(x)|. \qquad (3.3)$$
But, by definitions of $\tb(x)$ and $t_n(x)$, we have for all $x\in\cv_n^t$ :
\begin{eqnarray*}
& & \frac{1}{\sqrt {nh^k}}|\tb(x)-t_n(x)|\\
 & \leq & \Bigg(|t-f(x)|\Bigg|\frac{1}{\sqrt {\kt f(x)}}-\frac{1}{\sqrt {V_n(x)h^{-k}}}\Bigg|+\sqrt {\frac{h^k}{V_n(x)}}|Ef_n(x)-f(x)|\Bigg)\\
& \leq & \Bigg( \frac{(\log n)^{\beta}}{\sqrt {nh^k}} \sqrt {\frac{|\kt f(x)-V_n(x)h^{-k}|}{\kt f(x)V_n(x)h^{-k}}}+\sqrt {\frac{h^k}{V_n(x)}}|Ef_n(x)-f(x)|\Bigg). \ (3.4)
\end{eqnarray*}
It is a classical exercise to deduce from {\bf H1}, {\bf H3} that, since $\cv_n^t$ is contained in $\cal C$, 
$$\sup_{x\in\cv_n^t} |Ef_n(x)-f(x)|\leq c\, h^2,$$
and similarly, that
$$\sup_{x\in\cv_n^t}|\kt f(x)-V_n(x)h^{-k}| \leq c\, h.$$
One deduces from (3.4) and above that
$$\sup_{x\in\cv_n^t} |\tb(x)-t_n(x)|\leq c\,\big(\sqrt h\, (\log n)^{\beta}+\sqrt {nh^{k+4}}\big).$$
Thus, by (3.3) and since $t\in\Theta_0(g)$, one has for all $n$ large enough :
$$E_n\leq c\,(\log n)^{\beta}\big(\sqrt h\,(\log n)^{\beta}+\sqrt {nh^{k+4}}\big),$$and the latter term vanishes by assumptions on $h$, hence the lemma $\bullet$

\bigskip
\noindent {\bf Proof of Proposition 3.1.} By Lemma 3.2, one only needs to prove Proposition 3.1 for all $t\in\Theta_0(g)$. Fix $t\in\Theta_0(g)$, and let
$$I_n:=\int_{\cv_n^t} \fb(t_n(x))d\lambda_g(x) \ {\rm and} \ {\overline I}_n:=\int_{\ov_n^t} \Phi(t_n(x))d\lambda_g(x).$$
By Lemma 3.3, the task is now to prove that
$$\lim_n \sqrt {nh^k}\, I_n=\sqrt {\frac{t\kt}{2\pi}} \int_{\partial \l(t)}\frac{g}{\|\nabla f\|} d\h=\lim_n \sqrt {nh^k}\, {\overline I}_n.$$
We only show the first equality. One has
$$I_n=\frac{1}{\sqrt {2\pi\kt}} \int_{\cv_n^t} 
\int_{b_n(x)}^{\infty} \exp\Big(-\frac {u^2}{2\kt}\Big)du \, d\lambda_g(x),$$
where for all $x\in\R^k$ such that $f(x)>0$, $b_n(x)=\sqrt {nh^k}(t-f(x))/f(x)^{1/2}$. By Fubini's theorem :
$$I_n=\frac{1}{\sqrt {2\pi\kt}}\int_0^{\infty} \exp\Big(-\frac{u^2}{2\kt}\Big)\lambda_g\Big(f^{-1}\Big[\max\Big(t-\frac{(\log n)^{\beta}}{\sqrt {nh^k}},\chi\Big(\frac{u}{\sqrt {nh^k}}\Big)^2\Big),t\Big]\Big)du,$$
where for all $v\geq 0$, $\chi(v)=-v/2+(1/2)\sqrt {v^2+4t}$. It is straightforward to prove the equivalence :
$$u\in [0,r_n] \Leftrightarrow \chi\Big(\frac{u}{\sqrt {nh^k}}\Big)^2\geq t-\frac{(\log n)^{\beta}}{\sqrt {nh^k}},$$
where $r_n=(\log n)^{\beta}/\sqrt {t-(\log n)^{\beta}(nh^k)^{-1/2}}$, so that one can split $I_n$ into two terms, {\it i.e.}, $I_n=I_n^1+I_n^2$, where 
\begin{eqnarray*}
I_n^1 & = & \frac{1}{\sqrt {2\pi\kt}} \int_0^{r_n} \exp\Big(-\frac{u^2}{2\kt}\Big) \lambda_g\Big(f^{-1}\Big[\chi\Big(\frac{u}{\sqrt {nh^k}}\Big)^2,t\Big]\Big)du \\
{\rm and} \ I_n^2 & = & \frac{1}{\sqrt {2\pi\kt}} \int_{r_n}^{\infty} \exp\Big(-\frac{u^2}{2\kt}\Big) \lambda_g\Big( f^{-1}\Big[t-\frac{(\log n)^{\beta}}{\sqrt {nh^k}},t\Big]\Big)du.
\end{eqnarray*}
Since $t\in\Theta_0(g)$, one has for all $n$ large enough :
$$\sqrt {nh^k}\, I_n^2\leq c\,(\log n)^{\beta} \int_{r_n}^{\infty} \exp\Big(-\frac{u^2}{2\kt}\Big)du, \qquad (3.5)$$
and the rightmost term vanishes. Thus, it remains to compute the limit of $\sqrt {nh^k} I_n^1$. Using an expansion of $\chi$ in a neighborhood of the origin, we get 
$$\lim_n \sqrt {nh^k}\, \lambda_g\Big(f^{-1} \Big[\chi\Big(\frac{u}{\sqrt {nh^k}}\Big)^2,t\Big]\Big)=u\sqrt t \int_{\partial \l(t)} \frac{g}{\|\nabla f\|} d\h, \qquad (3.6)$$
for all $u\geq 0$, since $t\in\Theta_0(g)$. Moreover, one deduces from Lemma 3.2 that for all $n$ large enough and for all $u\in [0,r_n]$ :
\begin{eqnarray*}
\sqrt {nh^k}\, \lambda_g\Big(
f^{-1}\Big[\chi\Big(\frac{u}{\sqrt {nh^k}}\Big)^2,t\Big]\Big) & \leq &  c\sqrt {nh^k} \Big(t-\chi\Big(\frac{u}{\sqrt {nh^k}}\Big)^2\Big)\\
& \leq & c\, u, \qquad (3.7)
\end{eqnarray*}
because $r_n/\sqrt {nh^k}\to 0$. Thus, according to (3.5)-(3.7) and the Lebesgue theorem :
\begin{eqnarray*}
\lim_n \sqrt {nh^k}\, I_n & = & \lim_n \sqrt {nh^k}\, I_n^1 \\
& = & \frac{1}{\sqrt {2\pi\kt}} \int_0^{\infty} 
\exp \Big(-\frac{u^2}{2\kt}\Big) u \sqrt t \int_{\partial \l(t)} 
\frac{g}{\|\nabla f\|}d\h\, du\\
& = & \sqrt {\frac{t\kt}{2\pi}} \int_{\partial \l(t)} \frac{g}{\|\nabla f\|}d\h,
\end{eqnarray*}
hence the proposition $\bullet$

\bigskip
\noindent {\bf 3.3.  Proof of Proposition 3.2.} From now on, we introduce two random variables $N_1$, $N_2$ with law ${\cal N} (0,1)$ such that $N_1,N_2,X_1,X_2,\cdots$ are independent. We let 
$$\sigma_n=\frac{1}{(\log n)^{2\beta}\log \log n}, \ \forall n\geq 2.$$
(As we will see later, the random variable $Z_n(x)+\sigma_n N_1$ -for instance- has a density with respect to the Lebesgue measure.) For simplicity, we assume in the following that under {\bf H3}, the support of $K$ is contained in the euclidean unit ball of $\R^k$. 

\bigskip
\noindent {\bf Lemma 3.4.} {\it Let $g\, :\, \R^k\to\R_+$ be a borel bounded function and assume that {\bf H2} holds. If $nh^k/(\log n)^{2\beta}\to\infty$, then for all $t\in\Theta_0(g)$ there exists $c>0$ such that for $n$ large enough :
\begin{eqnarray*}
& & \int_{\cv_n^t} P\Big(\Big\{ Z_n(x)\geq \tb(x)\Big\}\Delta\Big\{Z_n(x)+\sigma_n N_1\geq \tb(x)\Big\}\Big) d\lambda_g(x)\leq c\, w_n; \\
& and & \int_{\ov_n^t} P\Big(\Big\{ Z_n(x)< \tb(x)\Big\}\Delta\Big\{Z_n(x)+\sigma_n N_1< \tb(x)\Big\}\Big) d\lambda_g(x)\leq c\, w_n,
\end{eqnarray*}
where $w_n=(\log n)^{\beta}/(nh^k)+\sigma_n (\log n)^{\beta}/\sqrt {nh^k}$.}

\bigskip
\noindent {\bf Proof.} We only prove the first inequality. Let  $t\in\Theta_0(g)$ and 
$$P_n := \int_{\cv_n^t} P\Big(\Big\{ Z_n(x)\geq \tb(x)\Big\}\Delta\Big\{Z_n(x)+\sigma_n N_1\geq \tb(x)\Big\}\Big) d\lambda_g(x).$$ By independence of $N_1$ and $Z_n(x)$, $P_n$ is smaller than 
$$ \int_{\cv_n^t} \int \exp\Big(-\frac{z^2}{2}\Big) P\Big(\Big\{ Z_n(x)\geq \tb(x)\Big\}\Delta\Big\{Z_n(x)+\sigma_n z\geq \tb(x)\Big\}\Big) dz\, d\lambda_g(x),$$
and consequently,
$$P_n\leq \int_{\cv_n^t} \int\exp\Big(-\frac{z^2}{2}\Big) P\Big(|Z_n(x)-\tb(x)|\leq \sigma_n |z|\Big) dz\, d\lambda_g(x).$$
Since $t\in\Theta_0(g)$, one deduces from Lemma 3.1 that for $n$ large enough :
\begin{eqnarray*}
P_n & \leq & c\, \frac{\lambda_g(\cv_n^t)}{\sqrt {nh^k}}+\int_{\cv_n^t} \int \exp\Big(-\frac{z^2}{2}\Big) P\Big(|N_1-\tb(x)|\leq \sigma_n|z|\Big) dz\, d\lambda_g(x)\\
& \leq & c\,\Big( \frac{(\log n)^{\beta}}{nh^k}+\frac{\sigma_n (\log n)^{\beta}}{\sqrt {nh^k}}\Big),
\end{eqnarray*}
hence the lemma $\bullet$

\bigskip
\noindent {\bf Lemma 3.5.} {\it Fix $t\in\Theta$ and assume that {\bf H1}, {\bf H3} hold. Then, there exists a polynomial function $Q$ of degree 5 defined on $\R^2$  such that for all $(u_1,u_2)\in\R^2$ and $n$ large enough :
$$\Big|E\exp\Big(i\Big(u_1 Z_n(x)+u_2Z_n(y)\Big)\Big)-E\exp\Big(iu_1Z_n(x)\Big)E\exp\Big(iu_2Z_n(y)\Big)\Big|$$
$$\leq \frac{Q(|u_1|,|u_2|)}{\sqrt {nh^k}},$$
if $x,y\in\cv_n^t\cup\ov_n^t$ are such that $\|x-y\|\geq 2h$.
}

\bigskip
\noindent {\bf Proof.} First of all, fix $u_1,u_2\in\R$, $x,y\in\cv_n^t\cup\ov_n^t$ and consider the following quantities :
\begin{eqnarray*}
M_1 & := & \frac{u_1}{\sqrt {nV_n(x)}}\Big[K\Big(\frac{x-X}{h}\Big)-E K\Big(\frac{x-X}{h}\Big)\Big]\\
 {\rm and} \ M_2 & := & \frac{u_2}{\sqrt {nV_n(y)}}\Big[K\Big(\frac{y-X}{h}\Big)-E K\Big(\frac{y-X}{h}\Big)\Big].
\end{eqnarray*}
One deduces from the inequality $|\exp(iw)-1-iw+w^2/2|\leq |w|$ $\forall w\in\R$ that
$$\Big|E\exp\Big(i\Big(M_1+M_2\Big)\Big)-1+\frac{1}{2}E(M_1+M_2)^2\Big|$$
$$= \Big|E\Big[\exp\Big(i\Big(M_1+M_2\Big)\Big)-1-i(M_1+M_2)+\frac{1}{2}(M_1+M_2)^2\Big]\Big|\leq E|M_1+M_2|^3.$$
In a similar fashion, if $j=1$ or $j=2$ :
$$\Big|E\exp(iM_j)-1+\frac{1}{2}EM_j^2\Big|=\Big|E\Big[\exp(iM_j)-1-iM_j+\frac{1}{2} M_j^2\Big]\Big|\leq E|M_j|^3.$$
Consequently,
\begin{eqnarray*}
& & \Big|E\exp\Big(i\Big(M_1+M_2\Big)\Big)-E\exp\Big(iM_1\Big)E\exp\Big(iM_2\Big)\Big|\\
& \leq & E|M_1+M_2|^3+\Big|\Big(1-\frac{1}{2}E|M_1+M_2|^2\Big)-\Big(1-\frac{1}{2}EM_1^2\Big)\Big(1-\frac{1}{2}EM_2^2\Big)\Big|\\
& & + \Big|1-\frac{1}{2}EM_1^2\Big|E|M_2|^3+ \Big|1-\frac{1}{2}EM_2^2\Big|E|M_1|^3. \qquad (3.8)
\end{eqnarray*}
It is an easy exercice to prove that for all $n$ large enough, one has $\inf V_n(x)\geq ch^k$, the infinimum being taken over all $x\in\cv_n^t\cup\ov_n^t$. Consequently, if $j=1$ or $j=2$ :
$$E|M_j|^3\leq c\,\frac{|u_j|^3}{\sqrt {n^3h^k}},$$
from which we deduce that :
$$E|M_1+M_2|^3\leq c\,\frac{|u_1|^3+|u_2|^3}{\sqrt {n^3h^k}}.$$
Moreover, $EM_1^2=u_1^2/n$, $EM_2^2=u_2^2/n$ and for all $x,y\in\cv_n^t\cup\ov_n^t$ such that $\|x-y\|\geq 2h$ :
$$E(M_1+M_2)^2=EM_1^2+EM_2^2-\frac{u_1u_2}{n\sqrt {V_n(x)V_n(y)}} EK\Big(\frac{x-X}{h}\Big)EK\Big(\frac{y-X}{h}\Big),$$
because the support of $K$ is contained in the unit ball and hence
$$EK\Big(\frac{x-X}{h}\Big)K\Big(\frac{y-X}{h}\Big)=0.$$
 One deduces from above and (3.8) that for all $x,y\in\cv_n^t\cup\ov_n^t$ such that $\|x-y\|\geq 2h$ :
\begin{eqnarray*}
& & \Big|E\exp\Big(i\Big(M_1+M_2\Big)\Big)-E\exp\Big(iM_1\Big)E\exp\Big(iM_2\Big)\Big|\\
& \leq & c\,\frac{|u_1|^3+|u_2|^3}{\sqrt {n^3h^k}}+\frac{(u_1u_2)^2}{n^2}+c\,\frac{|u_2|^3(1+u_1^2)+|u_1|^3(1+u_2^2)}{\sqrt {n^3h^k}}+c\,\frac{|u_1u_2|h^k}{n}.
\end{eqnarray*}
By assumption, $nh^{3k}\to 0$ so that for $n$ large enough : $h^k\leq 1/\sqrt {nh^k}$. Consequently,
$$\Big|E\exp\Big(i\Big(M_1+M_2\Big)\Big)-E\exp\Big(iM_1\Big)E\exp\Big(iM_2\Big)\Big|\leq \frac{Q(|u_1|,|u_2|)}{\sqrt {nh^k}},$$
where $Q$ is defined for all $u_1,u_2\in\R$ by :
$$Q(u_1,u_2)=c\big(u_1^3+u_2^3+(u_1u_2)^2+u_1u_2+u_2^2u_1^3+u_1^3u_2^2\big).$$
Consequently, for all $u_1,u_2\in\R$ and $x,y\in\cv_n^t\cup\ov_n^t$ such that $\|x-y\|\geq 2h$ :
\begin{eqnarray*}
& & \Big|E\exp\Big(i\Big(u_1Z_n(x)+u_2Z_n(y)\Big)\Big)-E\exp\Big(iu_1Z_n(x)\Big)E\exp\Big(iu_2Z_n(y)\Big)\Big|\\
& = & \Big|\Big(E\exp\Big(i\Big(M_1+M_2\Big)\Big)\Big)^n-\Big(E\exp\Big(iM_1\Big)E\exp\Big(iM_2\Big)\Big)^n\Big|\\
& \leq & n\Big|E\exp\Big(i\Big(M_1+M_2\Big)\Big)-E\exp\Big(iM_1\Big)E\exp\Big(iM_2\Big)\Big|\\
& \leq & \frac{Q(|u_1|,|u_2|)}{\sqrt {nh^k}},
\end{eqnarray*}
 hence the lemma $\bullet$

\bigskip
 In the following, $uv$ stands for the usual scalar product of $u,v\in\R^2$.

\bigskip
\noindent {\bf Lemma 3.6.} {\it Let $x,y\in\R^k$ be such that $V_n(x)V_n(y)\neq 0$. Then, the bivariate random variable 
\[ \left( \begin{array}{clcr} Z_n(x)+\sigma_n N_1 \\ Z_n(y)+\sigma_n N_2 \\ \end{array} \right) \]
has a density $\varphi_n^{x,y}$ defined for all $u\in\R^2$ by
$$\varphi_n^{x,y}(u)=\frac{1}{4\pi^2} \int E\Big[\exp\Big(i\Big(v_1Z_n(x)+v_2Z_n(y)\Big)\Big)\Big]\exp\Big(-i\, uv-\frac{1}{2}\sigma_n^2\|v\|^2\Big)dv.$$
}

\bigskip
\noindent {\bf Proof.} By independence of $X_1,\cdots,X_n,N_1$ and $N_2$, the random variable 
\[ \left( \begin{array}{clcr} Z_n(x)\\ Z_n(y)\\ \end{array} \right) +
\sigma_n \left( \begin{array}{clcr} N_1\\ N_2\\ \end{array} \right) \]
has a density  $\varphi_n^{x,y}$ defined for all $u=(u_1,u_2)\in\R^2$ by 
$$\varphi_n^{x,y}(u)  = \frac{1}{2\pi\sigma_n^2} E\Big[\exp\Big(-\frac{(u_1-Z_n(x))^2}{2\sigma_n^2}\Big)\exp\Big(-\frac{(u_2-Z_n(y))^2}{2\sigma_n^2}\Big)\Big].$$
Using the equality
$$\frac{1}{\sqrt {2\pi\sigma_n^2}}\exp\Big(-\frac{z^2}{2\sigma_n^2}\Big)=\frac{1}{2\pi} \int \exp\Big(-izw-\frac{1}{2}\sigma_n^2w^2\Big)dw \ \forall z\in\R,$$
we deduce from the Fubini theorem that 
$$\varphi_n^{x,y}(u) = \frac{1}{4\pi^2} \int E\Big[\exp\Big(i\Big(v_1Z_n(x)+v_2Z_n(y)\Big)\Big)\Big] \exp\Big(-iuv-\frac{1}{2}\sigma_n^2\|v\|^2\Big)dv,$$
hence the lemma $\bullet$

\bigskip
\noindent {\bf Proof of Proposition 3.2.} We only prove the first equality of Proposition 3.2. According to Lemma 3.2, one only needs to prove the result for each $t\in\Theta_0(g)$. Hence we fix $t\in\Theta_0(g)$ and we put :
$$A_n(x)=\Big\{Z_n(x)\geq \tb(x)\Big\}, \  A_n^j(x)=\Big\{Z_n(x)+\sigma_nN_j\geq \tb(x)\Big\}, \ j=1,2,$$
for all $x\in\R^k$ such that $V_n(x)\neq 0$. First note that since the events $A_n(x)$ and $\{f_n(x)\geq t\}$ are equal, one has
\begin{eqnarray*}
&  &\var\Big[\lambda_g\Big(\cv_n^t\cap\l_n(t)\Big)\Big]\\
& = & \int_{(\cv_n^t)^{\times 2}} \Big(P(A_n(x)\cap A_n(y))-P(A_n(x))P(A_n(y))\Big)
d\lambda_g^{\otimes 2}(x,y). \quad (3.9)
\end{eqnarray*}
But, by Lemma 3.4 and since $t\in\Theta_0(g)$, one has for all $n$ large enough :
\begin{eqnarray*}
& & nh^k \int_{(\cv_n^t)^{\times 2}} \Big(P(A_n(x)\cap A_n(y))-P(A_n^1(x)\cap A_n^2(y))\Big)d\lambda_g^{\otimes 2}(x,y)\\
& \leq & 2nh^k \lambda_g(\cv_n^t)\int_{\cv_n^t} P(A_n(x)\Delta A_n^1(x)) d\lambda_g(x)\\
& \leq & c\,(\log n)^{\beta} \sqrt {nh^k} \Big( \frac{(\log n)^{\beta}}{nh^k}+\frac{\sigma_n(\log n)^{\beta}}{\sqrt {nh^k}}\Big)\\
& \leq & c\,\Big( \frac{(\log n)^{2\beta}}{\sqrt {nh^k}}+\sigma_n (\log n)^{2\beta}\Big),
\end{eqnarray*}
and the latter term tends to 0 by assumption. In a similar fashion, one can prove that 
$$nh^k \int_{(\cv_n^t)^{\times 2}} \Big(P(A_n(x))P(A_n(y))-P(A_n^1(x))P(A_n^2(y))\Big)d\lambda_g^{\otimes 2}(x,y)\to 0.$$
By the above results and (3.9), it remains to show that
$$nh^k\int_{(\cv_n^t)^{\times 2}} \Big( P(A_n^1(x)\cap A_n^2(y))-P(A_n^1(x))P(A_n^2(y)) \Big)d\lambda_g^{\otimes 2}(x,y) \to 0.\qquad (3.10)$$
Let $T(h)=\{(x,y)\in(\R^k)^{\times 2} :\, \|x-y\|\leq 2h\}$. According to the Fubini theorem,
\begin{eqnarray*}
nh^k \lambda_g^{\otimes 2}\Big((\cv_n^t)^{\times 2}\cap T(h)\Big) & = & nh^k \int_{\cv_n^t} \lambda_g\Big(\cv_n^t\cap B(x,2h)\Big)d\lambda_g(x)\\
& \leq & nh^k \int_{\cv_n^t} \lambda_g(B(x,2h))d\lambda_g(x),
\end{eqnarray*}
where $B(z,r)$ stands for the euclidean closed ball with center at $z\in\R^k$ and radius $r>0$. Since $t\in\Theta_0(g)$, one deduces that
\begin{eqnarray*}
nh^k \lambda_g^{\otimes 2}\Big((\cv_n^t)^{\times 2}\cap T(h)\Big) & \leq & c\, n h^k \frac{(\log n)^{\beta}}{\sqrt {nh^k}} h^k\\
& \leq & c\,\sqrt {nh^{3k} (\log n)^{2\beta}},
\end{eqnarray*}
so that, by assumption on the bandwidth $h$ :
$$\lim_n nh^k \lambda_g^{\otimes 2}\Big((\cv_n^t)^{\times 2}\cap T(h)\Big)=0.$$
Let now ${\cal S}_n=(\cv_n^t)^{\times 2}\cap T(h)^c$. According to (3.10) and the above result, one only needs now to prove that :
$$nh^k\int_{{\cal S}_n}\Big(P(A_n^1(x)\cap A_n^2(y))-P(A_n^1(x))P(A_n^2(y))\Big)d\lambda_g^{\otimes 2}(x,y)\to 0.\quad (3.11)$$
By Lemmas 3.5 and 3.6, one has for all $x,y\in{\cal S}_n$:
\begin{eqnarray*}
& & \Big|P(A_n^1(x)\cap A_n^2(y))-P(A_n^1(x))P(A_n^2(y))\Big| \\
& \leq & \int \Big|E\exp\Big(i\Big(u_1 Z_n(x)+u_2Z_n(y)\Big)\Big)\\
& & \qquad -E\exp\Big(iu_1Z_n(x)\Big)E\exp\Big(iu_2Z_n(y)\Big)\Big|
\exp\Big(-\frac{1}{2}\sigma_n^2\|u\|^2\Big)du_1du_2\\
& \leq & \frac{1}{\sqrt {nh^k}}\int Q(|u_1|,|u_2|)\exp\Big(-\frac{1}{2}\sigma_n^2\|u\|^2\Big)du_1du_2\\
& \leq & \frac{c}{\sigma_n^7 \sqrt {nh^k}},
\end{eqnarray*}
where $Q$ is the polynomial function defined in Lemma 3.5. Consequently, one has for all $n$ large enough :
\begin{eqnarray*}
& & nh^k\int_{{\cal S}_n}\Big(P(A_n^1(x)\cap A_n^2(y))-P(A_n^1(x))P(A_n^2(y))\Big)d\lambda_g^{\otimes 2}(x,y)\\
& \leq & c\, \frac{\sqrt {nh^k}}{\sigma_n^7} \lambda_g^{\otimes 2}({\cal S}_n)\\
& \leq & c\,\frac{\sqrt {nh^k}}{\sigma_n^7} \lambda_g(\cv_n^t)^2\\
& \leq & c\,\frac{(\log n)^{2\beta}}{\sigma_n^7\sqrt {nh^k}},
\end{eqnarray*}
which tends to 0 by assumption, hence (3.11) $\bullet$

\bigskip
\noindent {\bf 4. Proof of Corollary 2.1.}

\bigskip
\noindent {\bf Lemma 4.1.} {\it  Let $k\geq 2$ and assume that {\bf H1}-{\bf H3} hold. If $nh^{k+4}(\log n)^2\to 0$ and  $nh^k/(\log n)^{16}\to\infty$, then for a.e. $t\in\Theta$ :
$$\sqrt {nh^k} \Big( \lambda_{f_n}(\l(t))-\lambda_{f_n}(\l_n(t))\Big)\lp 0.$$
}

\bigskip
\noindent {\bf Proof.} Let $t\in\Theta$ be such that the conclusion of Theorem 2.1 holds both for $g\equiv f$ and $g\equiv 1$. Notice that
\begin{eqnarray*}
\lambda_{f_n}(\l(t))-\lambda_{f_n}(\l_n(t)) & = & \int f_n\Big({\bf 1}_{\{f\geq t\}}-{\bf 1}_{\{f_n\geq t\}}\Big)d\lambda\\
& = & \int_{\l(t)} f_n {\bf 1}_{\{f_n<t\}}d\lambda-\int_{\l(t)^{c}} f_n {\bf 1}_{\{f_n\geq t\}}d\lambda.
\end{eqnarray*}
As in the proof of Theorem 2.1, we see that the result of the lemma will hold if we show that $\sqrt {nh^k} K_n\lp 0$, where
$$K_n := \int_{\ov_n^t} f_n{\bf 1}_{\{f_n<t\}}d\lambda-\int_{\cv_n^t} f_n{\bf 1}_{\{f_n\geq t\}}d\lambda.$$
Split $K_n$ into four terms as follows :
\begin{eqnarray*}
K_n & = &\int_{\ov_n^t} (f_n-f){\bf 1}_{\{f_n<t\}}d\lambda-\int_{\cv_n^t} (f_n-f){\bf 1}_{\{f_n\geq t\}}d\lambda\\
& & +\int_{\ov_n^t}{\bf 1}_{\{f_n<t\}}d\lambda_f-\int_{\cv_n^t}{\bf 1}_{\{f_n\geq t\}}d\lambda_f. \qquad (4.1)
\end{eqnarray*}
On one hand, it is a classical exercise to deduce from {\bf H1}, {\bf H3} that 
$$\sup_{\ov_n^t} |f_n-f|\lp 0.$$
 Thus, using (3.2),
$$\sqrt {nh^k} \int_{\ov_n^t} (f_n-f){\bf 1}_{\{f_n<t\}}d\lambda \lp 0.$$
 In a similar fashion :
$$\sqrt {nh^k} \int_{\cv_n^t} (f_n-f){\bf 1}_{\{f_n\geq t\}}d\lambda \lp 0.$$
On the other hand, we get from (3.2) that :
$$\lim_n \sqrt {nh^k} \int_{\cv_n^t} {\bf 1}_{\{f_n\geq t\}}d\lambda_f=\lim_n \sqrt {nh^k} \int_{\ov_n^t} {\bf 1}_{\{f_n < t\}}d\lambda_f,$$
where the limits are in probability. By  the above results and (4.1), 
$\sqrt {nh^k}K_n$ tends to 0 in probability, hence the lemma $\bullet$

\bigskip
\noindent {\bf Lemma 4.2.} {\it  Let $k\geq 2$, $t\in\Theta$ and assume that {\bf H1}, {\bf H3} hold. If $nh^{k+4}\to 0$, then :
$$\sqrt {nh^k} \Big( \lambda_f(\l(t))-\lambda_{f_n}(\l(t))\Big)\lp 0.$$
}

\bigskip
\noindent {\bf Proof.} Observe that 
$$\lambda_f(\l(t))-\lambda_{f_n}(\l(t))=\int_{\l(t)} (f-Ef_n)d\lambda+\int_{\l(t)}(Ef_n-f_n)d\lambda.$$
According to {\bf H1}, {\bf H3}, we have :
$$\int_{\l(t)} |f-Ef_n|d\lambda\leq ch^2,$$
and since $nh^{k+4}\to 0$, we only need to prove that 
$$\sqrt {nh^k} \int_{\l(t)} (Ef_n-f_n)d\lambda\lp 0.$$
We prove that this convergence holds in quadratic mean. We have :
\begin{eqnarray*}
E\Big(\sqrt {nh^k} \int_{\l(t)} (Ef_n-f_n)d\lambda\Big)^2
 & \leq & \frac{1}{h^k} E\Big( \int_{\l(t)} K\Big(\frac{x-X}{h}\Big)dx\Big)^2\\
& \leq & \frac{1}{h^k} \int_{\l(t)^{\times 2}} EK\Big(\frac{x-X}{h}\Big)K\Big(\frac{y-X}{h}\Big) dxdy.
\end{eqnarray*}
Recall that we assume in Section 3.3 that the support of $K$ is contained in the unit ball so that if $\|x-y\|\geq 2h$,
$$EK\Big(\frac{x-X}{h}\Big)K\Big(\frac{y-X}{h}\Big)=0.$$
Letting  $R(h)=\{(x,y)\in\l(t)^{\times 2} :\, \|x-y\|\leq 2h\}$, one deduces from above that
\begin{eqnarray*}
E\Big(\sqrt {nh^k} \int_{\l(t)} (Ef_n-f_n)d\lambda\Big)^2 & \leq & 
\frac{c}{h^k} \int_{R(h)} \int K\Big(\frac{x-u}{h}\Big)f(u)dudxdy\\
& \leq & c\int_{R(h)}  \int K(v)f(x-hv)dvdxdy\\
& \leq & c\, \lambda^{\otimes 2} (R(h))\\
& \leq & c \int_{\l(t)} \lambda\Big(\l(t)\cap B(x,2h)\Big)dx,
\end{eqnarray*}
according to the Fubini theorem. Thus, we get :
$$E\Big(\sqrt {nh^k} \int_{\l(t)} (Ef_n-f_n)d\lambda\Big)^2 \leq ch^k,$$
hence the lemma $\bullet$

\bigskip
\noindent {\bf Lemma 4.3.} {\it Let $p\in [0,1]$ and assume that {\bf H1}, {\bf H3} and {\bf H4} hold. If $nh^k/\log n\to \infty$, then $t_n^{(p)}\to t^{(p)}$ a.s.}

\bigskip
\noindent {\bf Proof.} Let $t=t^{(p)}$ and $t_n=t_n^{(p)}$. As seen in the proof of Theorem 2.1, $\sup_{\R^k} |f_n-f|\to 0$ a.s. Hence, one can fix 
$$\omega\in\Big\{\sup_{\R^k}|f_n-f|\to 0\Big\}.$$
For notational convenience, we omit $\omega$ until the end of this proof. Since $f$ is bounded, one has $\sup_n\sup_{\R^k} f_n<\infty$ and consequently $\sup_n t_n<\infty$. Thus, from each sequence of integers, one can extract a subsequence $(n_k)_k$ such that $t_{n_k}\to t^*$. On one hand, according to Scheffé's theorem, 
$$\lim_n \Big(\lambda_{f_{n_k}}(\l_{n_k}(t_{n_k}))-\lambda_f(\l_{n_k}(t_{n_k}))\Big)=0, \qquad (4.2)$$
since both $f$ and $f_{n_k}$ are density functions on $\R^k$ and 
$$\Big|\lambda_{f_{n_k}}(\l_{n_k}(t_{n_k}))-\lambda_{f}(\l_{n_k}(t_{n_k}))\Big|\leq \int |f_{n_k}-f|d\lambda.$$
On the other hand, letting $\eps_k=\sup_{\R^k}|f_{n_k}-f|$, one observes that
\begin{eqnarray*}
\Big|\lambda_f(\l(t_{n_k}))-\lambda_f(\l_{n_k}(t_{n_k}))\Big| & =  &
\int f\Big| {\bf 1}_{\{f\geq t_{n_k}\}}-{\bf 1}_{\{ f_{n_k}\geq t_{n_k}\}}\Big| d\lambda\\
& \leq & \int f {\bf 1}_{\{t_{n_k}-\eps_k\leq f\leq t_{n_k}+\eps_k\}}d\lambda\\
& \leq & c\,\lambda\Big( f^{-1}([t_{n_k}-\eps_k,t_{n_k}+\eps_k]\cap (0,\sup_{\R^k}f])\Big),
\end{eqnarray*}
and the latter term tends to 0 as $k\to\infty$  under {\bf H4} (consider separately the two cases : $t^*=0$ and $t^*>0$). One deduces from (4.2) that :
\begin{eqnarray*}
\lim_n \Big(\lambda_f(\l(t))-\lambda_f(\l(t_{n_k}))\Big) & = & \lim_n \Big(p-\lambda_f(\l(t_{n_k}))\Big) \\
& = & \lim_n \Big( \lambda_{f_{n_k}}(\l_{n_k}(t_{n_k}))-\lambda_f(\l_{n_k}(t_{n_k}))\Big)\\
& & + \lim_n \Big(\lambda_f(\l_{n_k}(t_{n_k}))-\lambda_f(\l(t_{n_k}))\Big)\\
& = & 0. \qquad (4.3)
\end{eqnarray*}
Moreover, the application $s\mapsto\lambda_f(\l(s))$ defined on $[0,\sup_{\R^k}f]$ is continuous according to {\bf H4}. Consequently, one has 
$$\lim_n \lambda_f(\l(t_{n_k}))=\lambda_f(\l(t^*)),$$
and thus, by (4.3), $\lambda_f(\l(t))=\lambda_f(\l(t^*))$ and hence $t=t^*$ because $\cal P$ is one-to-one. One conclude $t_n\to t$ since we proved that from each sequence of integers, one can extract a subsequence $(n_k)_k$ such that $t_{n_k}\to t$. The lemma is proved $\bullet$

\bigskip
\noindent {\bf Lemma 4.4.} {\it Let $k\geq 2$ and assume that {\bf H1}-{\bf H4} hold. If $nh^{k+4}(\log n)^2\to 0$ and $nh^{k+2}/\log n\to \infty$, then for a.e. $p\in {\cal P}(\Theta)$ :
$$\sqrt {nh^k} \int_{t_n^{(p)}}^{t^{(p)}} \int_{\partial \l_n(s)} \frac{1}{\|\nabla f_n\|}d\h\, ds\lp 0.$$
}

\bigskip
\noindent {\bf Proof.} One only needs to choose $p\in{\cal P}(\Theta)$ such that the conclusion of Lemma 4.1 holds for $t^{(p)}$. For simplicity, let $t=t^{(p)}$ and $t_n=t_n^{(p)}$. It is a classical exercise to prove that since $nh^{k+2}/\log n\to \infty$ and $nh^{k+4}\to 0$,
$$\|\nabla f_n\|\to \|\nabla f\| \ {\rm a.s.},$$
uniformly over the compact sets. Thus, by Lemma 4.3 and {\bf H2}, we have a.s. and for $n$ large enough :
$$\inf_{f^{-1}[\min(t_n,t),\max(t_n,t)]} \|\nabla f_n\|>0. \qquad (4.4)$$
We deduce from Proposition A that a.s. and for $n$ large enough :
\begin{eqnarray*}
\lambda_{f_n}(\l_n(t_n))-\lambda_{f_n}(\l_n(t)) & = & \int \Big( {\bf 1}_{\{f_n\geq t_n\}}-{\bf 1}_{\{f_n\geq t\}}\Big)d\lambda_{f_n}\\
& = & \int {\bf 1}_{\{t_n\leq f_n<t\}}d\lambda_{f_n}-\int {\bf 1}_{\{t\leq f_n < t_n\}}d\lambda_{f_n}\\
& = & \int_{t_n}^t \int_{\partial \l_n(s)} \frac{f_n}{\|\nabla f_n\|} d\h\, ds,
\end{eqnarray*}
where the latter integral is defined according to (4.4). Consequently, 
$$\Big|\lambda_{f_n}(\l_n(t_n))-\lambda_{f_n}(\l_n(t))\Big| = \int_{\min(t_n,t)}^{\max(t_n,t)} s \int_{\partial \l_n(s)} 
\frac{1}{\|\nabla f_n\|} d\h\, ds.$$
By Lemma 4.3, one has a.s. and for $n$ large enough : $t_n\geq t/2$. Since $\lambda_{f_n}(\l_n(t_n))=p=\lambda_f(\l(t))$, one deduces that :
$$\Big|\lambda_f(\l(t))-\lambda_{f_n}(\l_n(t))\Big|\geq \frac{t}{2} \int_{\min(t_n,t)}^{\max(t_n,t)}\int_{\partial 
\l_n(s)} \frac{1}{\|\nabla f_n\|} d\h\, ds.$$
We can now conclude the proof of the lemma because
$$\sqrt {nh^k} \Big|\lambda_f(\l(t))-\lambda_{f_n}(\l_n(t))\Big|\lp 0,$$
by Lemmas 4.1 and 4.2 $\bullet$

\bigskip
\noindent {\bf Lemma 4.5} {\it Assume that {\bf H1}-{\bf H3} hold. If $nh^k/(\log n)^2\to\infty$, then for a.e. $p\in{\cal P}(\Theta)$ :
$$\frac{\sqrt {nh^k}}{\log n}\, |t_n^{(p)}-t^{(p)}|\lp 0.$$}

\bigskip
\noindent {\bf Proof.} By {\bf H2} and the Lebesgue-Besicovitch theorem (Evans and Gariepy, [16], Theorem 1, Chapter I), we have for a.e. $p\in{\cal P}(\Theta)$ :
$$\frac{1}{\eps} \int_{t^{(p)}-\eps}^{t^{(p)}} \int_{\partial \l(s)} \frac{f}{\|\nabla f\|} d\h\, ds\to \int_{\partial \l(t^{(p)})} \frac{f}{\|\nabla f\|} d\h,$$
as $\eps\searrow 0$. Thus, one only needs to prove the lemma for $p\in{\cal P}(\Theta)$ such that the above result holds. For convenience, let $t=t^{(p)}$ and $t_n=t_n^{(p)}$. It suffices to show that
$$\frac{\sqrt {nh^k}}{\log n}\, |t_n^{(p)}-t^{(p)}|\lp 0$$
on the event $A_n$ defined by
$$A_n=\Big\{ \sup_{\l(t/2)}|f_n-f|\leq r_n\Big\},$$
where $r_n=(\log n)^{3/4}/\sqrt {nh^k}$, because $P(A_n)\to 1$ (see the proof of Theorem 2.1). According to Lemma 4.3, one has a.s. and for $n$ large enough : $\l(t_n)\cup\l_n(t_n)\subset \l(t/2)$ on the event $A_n$. Then, 
\begin{eqnarray*}
|\lambda_f(\l(t_n))-\lambda_{f_n}(\l_n(t_n))| & = & \Big| \int_{\l(t_n)} fd\lambda-\int_{\l_n(t_n)} f_nd\lambda\Big|\\
& \leq & \int_{\l(t/2)} |f_n-f|d\lambda+\int f\Big|{\bf 1}_{\l(t_n)}-{\bf 1}_{\l_n(t_n)}\Big|d\lambda\\
& \leq & c\, r_n+c\,\lambda\Big(\l(t_n)\Delta\l_n(t_n)\Big). \qquad (4.5)
\end{eqnarray*}
But, on $A_n$ :
$$\lambda\Big(\l(t_n)\Delta\l_n(t_n)\Big)\leq \lambda\Big(\Big\{t_n-r_n\leq f\leq t_n+r_n\Big\}\Big).$$
By {\bf H1}, {\bf H2}, there exists a neighborhood $V$ of $t$ such that 
$$\inf_{f^{-1} (V)} \|\nabla f\|>0,$$
thus, by Lemma 4.3, one has a.s. and for $n$ large enough :
\begin{eqnarray*}
\lambda\Big(\l(t_n)\Delta\l_n(t_n)\Big) & \leq & \sup_{s\in V} \lambda\Big(\Big\{s-r_n\leq f\leq s+r_n\Big\}\Big)\\
& \leq & c\,r_n,
\end{eqnarray*}
where the latter inequality is a consequence of Proposition A. According to (4.5), one has on $A_n$ and for $n$ large enough :
$$|\lambda_f(\l(t_n))-\lambda_f(\l(t))|=|\lambda_f(\l(t_n))-\lambda_{f_n}(\l_n(t_n))|\leq c\,r_n.$$
Observe now that by Proposition A and our choice of $t$, one has a.s. :
$$\frac{\lambda_f(\l(t_n))-\lambda_f(\l(t))}{t_n-t} \to \int_{\partial \l(t)} \frac{f}{\|\nabla f\|} d\h\neq 0,$$
thus on $A_n$, 
$$|t_n-t|\leq c\, r_n,$$
for $n$ large enough, hence the lemma $\bullet$

\bigskip
\noindent {\bf Lemma 4.6.} {\it Assume that {\bf H1}-{\bf H4} hold and let $(\alpha_n)_n$ be a sequence of positive real numbers. If $\alpha_n\to 0$, $\alpha_n^2nh^k/(\log n)^2\to\infty$ and $nh^k/(\log n)^2\to\infty$, then for a.e. $p\in{\cal P}(\Theta)$ :
$$\frac{1}{\alpha_n}\lambda\Big(\l_n(t_n^{(p)})-\l_n(t_n^{(p)}+\alpha_n)\Big)\lp \int_{\l(t^{(p)})} \frac{1}{\|\nabla f\|}d\h.$$}

\bigskip
\noindent {\bf Proof.} According to Proposition A and {\bf H1}, {\bf H2}, {\bf H4}, one has for a.e. $t\in\Theta$ :
$$\frac{1}{\eps}\lambda\Big(\l(t)-\l(t+\eps)\Big)=\frac{1}{\eps}\lambda\Big(\Big\{t\leq f\leq t+\eps\Big\}\Big)\to \int_{\l(t^{(p)})} \frac{1}{\|\nabla f\|}d\h,$$
as $\eps\searrow 0$. Hence, it suffices to prove the lemma for all $p\in{\cal P}(\Theta)$ such that the above result holds with $t=t^{(p)}$. For convenience, let $t=t^{(p)}$ and $t_n=t_n^{(p)}$. By Lemma 4.5, one only needs to prove that 
$$\frac{1}{\alpha_n}\lambda\Big(\l_n(t_n)-\l_n(t_n+\alpha_n)\Big)=\frac{1}{\alpha_n}\lambda\Big(\Big\{t_n\leq f_n < t_n+\alpha_n\Big\}\Big)\lp \int_{\l(t^{(p)})} \frac{1}{\|\nabla f\|}d\h,$$
on the event $B_n$ defined by
$$B_n=\Big\{\sup_{\l(t/2)}|f_n-f|\leq v_n, \ |t_n-t|\leq v_n\Big\},$$
where $v_n=\log n/\sqrt {nh^k}$, because $P(B_n)\to 1$. But, for $n$ large enough, one has $\l_n(t_n)\cup\l(t)\subset \l(t/2)$ on $B_n$. Consequently, 
$$\frac{1}{\alpha_n}\Big|\lambda\Big(\Big\{t_n\leq f_n<t_n+\alpha_n\Big\}\Big)-
\lambda\Big(\Big\{t\leq f\leq t+\alpha_n\Big\}\Big)\Big|$$
$$\leq \frac{1}{\alpha_n} \lambda\Big(\Big\{t-2v_n\leq f\leq t+2v_n\Big\}\Big) \leq c\,\frac{v_n}{\alpha_n},$$
and the latter term tends to 0 by assumption on $\alpha_n$. Finally, the choice of $t$ implies that 
$$\frac{1}{\alpha_n}\lambda\Big(\Big\{t\leq f_n\leq  t+\alpha_n\Big\}\Big)\to \int_{\l(t^{(p)})} \frac{1}{\|\nabla f\|}d\h,$$
so that on $B_n$ :
$$\frac{1}{\alpha_n}\lambda\Big(\Big\{t_n\leq f_n<t_n+\alpha_n\Big\}\Big)\lp \int_{\l(t^{(p)})} \frac{1}{\|\nabla f\|}d\h,$$
hence the lemma $\bullet$

\bigskip
\noindent {\bf Proof of Corollary 2.1.} According to Lemma 4.3, Lemma 4.6 and Theorem 2.1, one only needs to prove that for a.e. $p\in {\cal P}(\Theta)$ :
$$\sqrt {nh^k}\Big[\lambda\Big(\l_n(t_n^{(p)})\Delta \l(t^{(p)})\Big)-\lambda\Big(\l_n(t^{(p)})\Delta \l(t^{(p)})\Big)\Big]\lp 0.$$
  Moreover, it suffices to show the above result for each $p\in {\cal P}(\Theta)$ such that the conclusion of Lemma 4.4 holds. Fix such a $p\in{\cal P}(\Theta)$ and, for simplicity, let $t=t^{(p)}$ and $t_n=t_n^{(p)}$. A straightforward computation gives the relation :
$$D_n :=\lambda\Big(\l_n(t_n)\Delta \l(t)\Big)-\lambda\Big(\l_n(t)\Delta \l(t)\Big)=\int \Big({\bf 1}_{\{f_n\geq t_n\}}-{\bf 1}_{\{f_n\geq t\}}\Big)\eta\, d\lambda,$$
where $\eta=1-2{\bf 1}_{\{f\geq t\}}$. Then, 
$$D_n= \int {\bf 1}_{\{t_n\leq f_n<t\}}\eta\, d\lambda- \int {\bf 1}_{\{t\leq f_n< t_n\}}\eta\, d\lambda.$$
By (4.4) and {\bf H3}, one can now apply Proposition A, which gives :
$$D_n=\int_{t_n}^t \int_{\partial \l_n(s)} \frac{\eta}{\|\nabla f_n\|}d\h\, ds.$$
Consequently, 
$$|D_n|\leq \int_{\min(t_n,t)}^{\max (t_n,t)} \int_{\partial \l_n(s)} \frac{1}{\|\nabla f_n\|}d\h\, ds,$$
so that by Lemma 4.4 :
$$\sqrt {nh^k} D_n=\sqrt {nh^k} \Big[\lambda\Big(\l_n(t_n)\Delta \l(t)\Big)-\lambda\Big(\l_n(t)\Delta \l(t)\Big)\Big]\lp 0,$$
hence the corollary $\bullet$

\bigskip
\noindent {\bf Appendix : A change of variables formula.} Proposition A below is a consequence of the change of variables formula given in Evans and Gariepy ([16], Chapter III, Theorem 2). For a similar proof, see also Chapter III, Proposition 3 in the same book.

\bigskip
\noindent {\bf Proposition A.} {\it Let $\varphi\,:\, \R^k\to \R_+$ be a continuously differentiable function such that $\varphi(x)\to 0$ as $\|x\|\to\infty$, and  $I\subset \R_+$ be an interval such that $\inf I >0$ and
$$\inf_{\varphi^{-1}(I)} \|\nabla \varphi\|>0.$$
Then, for all borel bounded function  $g\, :\R^k\to \R$ :
$$\int_{\varphi^{-1}(I)} gdx=\int_I \int_{\varphi^{-1}(\{s\})} \frac{g}{\|\nabla \varphi\|} d\h\, ds.$$}

\bigskip
\noindent {\bf Proof.} Notice that $\varphi$ is a locally Lipschitz function and $$g{\bf 1}_{\varphi^{-1}(I)}$$
is integrable because $\varphi^{-1}(I)$ is bounded. Proposition A is then an easy consequence of Theorem 2 in Evans and Gariepy ([16], Chapter III) $\bullet$

\bigskip
\noindent {\bf Acknowledgements.} The author thank Andr\'e Mas and Nicolas Molinari for many helpful comments.

\bigskip
\bigskip
\centerline {\bf REFERENCES}

\bigskip
\bigskip
\smallskip
\noindent  [1] J.A. Hartigan, Clustering Algorithms, Wiley, New-York, 1975.

\smallskip
\noindent [2] L. Devroye and G.L. Wise, Detection of abnormal behavior via nonparametric estimation of the support, SIAM J. Appl. Math. 38 (1980) 480-488.

\smallskip
\noindent  [3] U. Grenander, Abstract Inference, Wiley, New-York, 1981.

\bigskip
\noindent  [4] A. Cuevas, On pattern analysis in the non-convex case, Kybernetes 19 (1990) 26-33.

\smallskip
\noindent  [5] A. Cuevas and R. Fraiman, Pattern analysis via nonparametric density estimation, Unpublished manuscript (1993).

\smallskip
\noindent  [6] I.S. Molchanov, Empirical estimation of distribution quantiles of random closed sets, Theory Probab. Appl. 35 (1990) 594-600.\smallskip

\noindent  [7] I.S. Molchanov, A limit theorem for solutions of inequalities, Unpublished manuscript (1993).

\smallskip
\noindent  [8] J.A. Hartigan, Estimation of a convex density contour in two dimensions, J. Amer. Statist. Assoc. 82 (1987) 267-270.

\smallskip
\noindent  [9] D.W. M\"uller, The excess mass approach in statistics, Beitr\"age zur Statistik, Univ. Heidelberg, 1993.

\smallskip
\noindent  [10] D.W. M\"uller and G. Sawitzki, Excess mass estimates and tests of multimodality, J. Amer. Statist. Assoc. 86 (1991) 738-746.

\smallskip
\noindent  [11] D. Nolan, The excess-mass ellipsoid, J. Multivariate Anal. 39 (1991) 348-371.

\smallskip
\noindent  [12] W. Polonik, Measuring mass concentration and estimating density contour clusters - an excess mass approach, Ann. Statist. 23 (1995) 855-881.

\smallskip
\noindent  [13] A.B. Tsybakov, On nonparametric estimation of density level sets, Ann. Statist. 25 (1997) 948-969.

\smallskip
\noindent  [14] M. Rosenblatt, Remarks on some nonparametric estimates of a density function, Ann. Math. Statist. 27 (1956) 832-837.

\smallskip
\noindent  [15] D. Pollard, Convergence of Stochastic Processes, Springer, New York, 1984.

\smallskip
\noindent  [16] L.C. Evans and R.F. Gariepy, Measure Theory and Fine Properties of Functions, CRC Press, Boca Raton, 1992.

\smallskip 
\noindent [17] W. Feller, An Introduction to Probability Theory and Its Applications, Wiley, New-York, 1992.

\end{document}